\documentclass[12pt,a4paper,leqno]{article}

\usepackage{amssymb,exscale}
\usepackage[centertags]{amsmath}
\usepackage{graphics}
\usepackage{theorem}

\numberwithin{equation}{section}

\theoremstyle{change}
\theorembodyfont{\rmfamily}
\newtheorem{theorem}[equation]{Theorem}
\newtheorem{lemma}[equation]{Lemma}

\newtheorem{Def}[equation]{Definition}

\newcommand{\sk}{\smallskip}

\renewcommand{\phi}{\varphi}

\renewcommand{\(}{\bigl(}
\renewcommand{\)}{\bigr)\vphantom{)}}

\newcommand{\bN}{\mathbf N}
\newcommand{\Om}{\Omega}
\newcommand{\om}{\omega}
\newcommand{\simp}{\mathrm{simp}}
\newcommand{\spi}{\mathrm{spider}}
\newcommand{\Var}{\mathrm{Var}\,}

\newcommand{\pd}{\partial}

\newcommand{\eps}{\varepsilon}

\newcommand{\De}{\Delta}

\newcommand{\dist}{{\,\mathrm{dist}\,}}

\newcommand{\rhomax}{\rho_{\mathrm{max}}}
\newcommand{\E}{\mathbb E}
\newcommand{\Prob}{\mathbb P}
\newcommand{\R}{\mathbb R}
\newcommand{\C}{\mathbb C}

\newcommand{\qed}{$\hfill\square$}

\begin{document}

\title{Trees, not cubes: hypercontractivity, cosiness, and noise stability}
\author{Oded Schramm \and Boris Tsirelson}
\date{}

\maketitle

\begin{abstract}
Noise sensitivity of functions on the leaves of a binary tree is studied,
and a hypercontractive inequality is obtained.  We deduce that the
spider walk is not noise stable.
\end{abstract}

\section*{Introduction}

\begin{figure}[b]
\begin{center}
\includegraphics{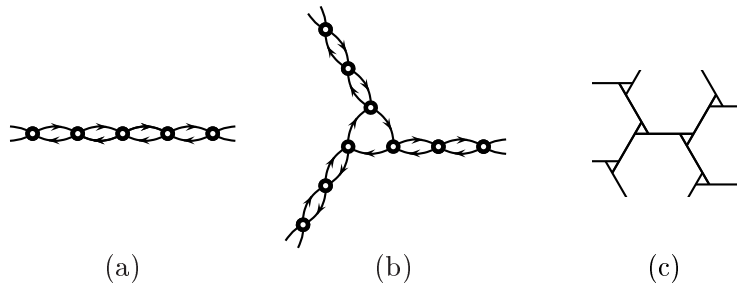}
\caption[]{(a) simple walk; (b) spider walk; (c) a spider web.
At each point, there are two equiprobable moves.}
\end{center}
\end{figure}

For the simplest random walk (Fig.~1a), the set $ \Om_n^\simp $ of all
$ n $-step trajectories may be thought of either as (the set of leaves
of) a binary tree, or (the vertices of) a binary cube $ \{-1,+1\}^n $.
However, consider
another random walk (Fig.~1b); call it the simplest \emph{spider}
walk, since it is a discrete counterpart of a spider martingale, see
\cite{BEKSY}. The corresponding $ \Om_n^\spi $ is the set of leaves of
a binary tree. For more complicated ``spider webs'' with several
``roundabouts'' we still have binary trees.  It is not quite
appropriate to think of such $n$-step ``spider walks'' as the vertices
of a binary cube, since for different $i$ and $j$ in $\{1,2,\dots,n\}$
it is not necessary that the $j$'th step has the same or opposite direction
from the $i$'th step.  Of course, one may choose to ignore this point,
and use the $n$ bits given by a point in $\{-1,1\}^n$ to describe a spider
walk, in such a way that for each $j=1,2,\dots,n$, the first
$j$ bits determine the first $j$ steps of the walk.  Such a correspondence
would not be unique.
In other words, cube structures on an $n$-level binary tree
may be introduced in different ways.

Noise sensitivity and stability are introduced and studied in
\cite{BKS} for functions on cubes. Different cube structures on a
binary tree are non-equivalent in that respect. It is shown here
that a natural function on $ \Om_n^\spi $ is non-stable
under \emph{every} cube structure.  One of the tools used is a new
hypercontractive inequality, which hopefully may find uses elsewhere.

\section{Stability and sensitivity on cubes, revisited}

A function $ f : \{-1,+1\}^n \to \C $ has its Fourier-Walsh expansion,
\begin{equation*}\begin{split}
& f (\tau_1,\dots,\tau_n) = \\
&\qquad = \hat f_0 + \sum_k \hat f_1 (k) \tau_k +
\sum_{k<l} \hat f_2 (k,l) \tau_k \tau_l + \dots + \hat f_n (1,\dots,n)
\tau_1 \dots \tau_n \, .
\end{split}\end{equation*}
Set 
\[
\tilde f_j(\tau_1,\dots,\tau_n) = \sum_{i_1<i_2<\cdots<i_j} \hat
f_j(i_1,\dots,i_j) \tau_{i_1}\tau_{i_2}\dots\tau_{i_j}\,.
\]
Since the transform $ f \mapsto \hat f $ is isometric, we have
$\| f \|^2 = \sum_{0}^n \| \tilde f_m \|^2$, where
\begin{equation}
\| f \|^2 = 2^{-n} \sum_{\tau_1,\dots,\tau_n} | f
  (\tau_1,\dots,\tau_n) |^2 \, . \label{norm} 
\end{equation}
The quantities
\[
S_1^m (f) = \sum_{i=1}^m \| \tilde f_i \|^2 \, , \quad S_m^\infty (f) =
\sum_{i=m}^n \| \tilde f_i \|^2 
\]
are used for describing low-frequency and high-frequency parts of the
spectrum of $ f $.

Given a sequence of functions $ F = \( f_n \)_{n=1}^\infty $, $ f_n :
\{-1,+1\}^n \to \C $, satisfying $ 0 < \liminf_{n\to\infty} \| f_n \|
\le \limsup_{n\to\infty} \| f_n \| < \infty $, we consider numbers
\[\begin{split}
S_1^m (F) &= \limsup_{n\to\infty} S_1^m (f_n) \, , \\
S_m^\infty (F) &= \limsup_{n\to\infty} S_m^\infty (f_n) \, .
\end{split}\]
Here is one of equivalent definitions of stability and sensitivity for
such $ F $, according to \cite[Th.~1.8]{BKS} (indicator functions are
considered there):
\begin{equation*}\begin{split}
F \text{ is stable} \quad &\text{iff} \quad S_m^\infty (F) \to 0
  \text{ for } m \to \infty \, , \\
F \text{ is sensitive} \quad &\text{iff} \quad S_1^m (F) = 0
  \text{ for all } m \, .
\end{split}\end{equation*}

A random variable $ \tau $ will be called a random sign, if $ \Prob
(\tau=-1) = 1/2 $ and $ \Prob (\tau=+1) = 1/2 $. A joint distribution
for two random signs $ \tau' $, $ \tau'' $ is determined by their
correlation coefficient $ \rho = \E ( \tau' \tau'' ) = 1 - 2 \Prob (
\tau' \ne \tau'' ) $. Given $ n $ independent pairs $ (\tau'_1,
\tau''_1), \dots, (\tau'_n, \tau''_n) $ of random signs with the same
correlation $ \rho $ for each pair, we call $ ( \tau'_1, \dots,
\tau'_n ) $ and $ ( \tau''_1, \dots, \tau''_n ) $ a \emph{$ \rho
$-correlated pair of random points of the cube $ \{-1,+1\}^n $.} (In
terms of \cite{BKS} it is $ \( x, N_\eps (x) \) $ with $ \eps =
(1-\rho)/2 $.) It is easy to see that
\[
\E \( \overline{f (\tau')} f (\tau'') \) = \sum_{m=0}^n \rho^m \|
\tilde f_m \|^2
\]
for a $ \rho $-correlated pair $ (\tau', \tau'') $. We may write it as
a scalar product in the space $ L_2 \( \{-1,+1\}^n \) $ with the norm
\eqref{norm},
\begin{equation}
\E \( \overline{f (\tau')} f (\tau'') \) = ( \rho^\bN f, f ) \, ;
\label{corr}
\end{equation}
here $ \rho^\bN $ is the operator $\rho^\bN f=\sum_n \rho^n\tilde f_n$.
Similarly, $ \E \( \overline{g(\tau')} f(\tau'') \) = ( \rho^\bN f, g
) $. On the other hand, 
\[ \E \( \overline{g(\tau')} f(\tau'') \) = \E \( \overline{g(\tau')}
\cdot \E \( f(\tau'') | \tau' \) \) = \( \tau' \mapsto \E \( f(\tau'')
| \tau' \), g \) \,;
\]
thus,
\begin{equation}
\E \( f(\tau'') | \tau' \) = ( \rho^\bN f) (\tau') \,
. \label{condexp}
\end{equation}
(Our $ \rho^\bN $ is $ T_\eta = Q_\eps $ of \cite{BKS} 
with $ \eta = \rho $, $ \eps = (1-\rho)/2 $.)
(In fact, let $\bN f =\sum_n n\tilde f_n$, then
 $ -\bN $ is the generator of a Markov process on $
\{-1,+1\}^n $; $ \exp (-t\bN) $ is its semigroup; note that $
\rho^\bN $ is of the form $ \exp (-t\bN) $. The Markov process is
quite simple: during $ dt $, each coordinate flips with the
probability $ \frac12 dt + o(dt) $. However, we do not need it.)
Note also that $ \E \( | f_n (\tau'') - (\rho^\bN f_n) (\tau') |^2 \,
\big| \, \tau' \) $ is the conditional variance $ \Var \( f_n (\tau'')
\big| \tau' \) $, and its mean value (over all $ \tau' $) is
\begin{equation}
\E \, \Var \( f_n (\tau'') \big| \tau' \) = \| f_n \|^2 - \| \rho^\bN f_n
\|^2 = \( (\mathbf1-\rho^{2\bN})f_n, f_n \) \, . \label{condvar}
\end{equation}
Note also that the operator $ 0^\bN = \lim_{\rho\to0} \rho^\bN $ is
the projection onto the one-dimensional space of constants, $ f
\mapsto (\E f) \cdot \mathbf 1 $.

Stability of $ F = \( f_n \)_{n=1}^\infty $ is equivalent to:
\begin{itemize}
\item $ \| \rho^\bN f_n - f_n \| \xrightarrow[\rho\to1]{} 0 $
uniformly in $ n $;
\item $ (\rho^\bN f_n, f_n) \xrightarrow[\rho\to1]{} \| f \|^2 $
uniformly in $ n $;
\item $ \| f_n \|^2 - \| \rho^\bN f_n \|^2 \xrightarrow[\rho\to1]{} 0 $
uniformly in $ n $.
\end{itemize}
Sensitivity of $ F $ is equivalent to:
\begin{itemize}
\item $ \| ( \rho^\bN - 0^\bN ) f_n \| \xrightarrow[n\to\infty]{} 0 $
for some (or every) $ \rho \in (0,1) $;
\item $ \( (\rho^\bN-0^\bN) f_n, f_n \) \xrightarrow[n\to\infty]{} 0 $
for some (or every) $ \rho \in (0,1) $.
\end{itemize}
Combining these facts with the probabilistic interpretation \eqref{corr},
\eqref{condexp}, \eqref{condvar} of $ \rho^\bN $ we see that
\begin{itemize}
\item $ F $ is stable iff $ \E \( \overline{f_n (\tau')} f_n (\tau'') \)
  \xrightarrow[\rho\to1]{} \E | f_n (\tau) |^2 $ uniformly in $ n $ or,
  equivalently, $ \E \( \Var ( f_n (\tau'') | \tau' )\)
  \xrightarrow[\rho\to1]{} 0 $ uniformly in $ n $;
\item $ F $ is sensitive iff $ \E \( \overline{f_n (\tau')} f_n (\tau'') \) -
  \big| \E f_n (\tau) \big|^2 \xrightarrow[n\to\infty]{} 0 $ for some
  (or every) $   \rho \in (0,1) $ or, equivalently,
  $ \E \, \big| \E ( f (\tau'') | \tau' ) - \E f
  \big|^2 \xrightarrow[n\to\infty]{} 0 $ for some (or every) $ \rho \in
  (0,1) $.
\end{itemize}
These are versions of definitions introduced in \cite[Sect.~1.1,
1.4]{BKS}.

\section{Stability and sensitivity on trees}

A branch of the $ n $-level binary tree can be written as a sequence
of sequences $ () $, $ (\tau_1) $, $ (\tau_1,\tau_2) $, $
(\tau_1,\tau_2,\tau_3) $, \dots, $ (\tau_1,\dots,\tau_n) $. Branches
correspond to leaves $ (\tau_1,\dots,\tau_n) \in \{-1,+1\}^n
$. Automorphisms of the tree can be described as maps $ A :
\{-1,+1\}^n \to \{-1,+1\}^n $ of the form
\[\begin{split}
& A (\tau_1,\dots,\tau_n) = \( a() \tau_1, a(\tau_1) \tau_2,
a(\tau_1,\tau_2) \tau_3, \dots, a(\tau_1,\dots,\tau_{n-1})
\tau_n \)
\end{split}\]
for arbitrary functions $ a : \cup_{m=1}^n \{-1,+1\}^{m-1} \to
\{-1,+1\} $. (Thus, the tree has $ 2^1 \cdot 2^2 \cdot 2^4 \cdot
\ldots \cdot 2^{2^{n-1}} = 2^{2^n-1} $ automorphisms, while the cube $
\{-1,+1\}^n $ has only $ 2^n n! $ automorphisms.)

Here is an example of a tree automorphism (far from being a cube
automorphism):
\[
(\tau_1,\dots,\tau_n) \mapsto \( \tau_1, \tau_1 \tau_2, \dots,
\tau_1\dots\tau_n \) \, .
\]
The function $ f_n (\tau_1,\dots,\tau_n) = \frac1{\sqrt n}
(\tau_1+\dots+\tau_n) $ satisfies $ S_1^1 (f_n) = 1 $, $ S_2^\infty
(f_n) = 0 $. However, the function $ g_n (\tau_1,\dots,\tau_n) =
\frac1{\sqrt n} \( \tau_1 + \tau_1 \tau_2 + \dots + \tau_1\dots\tau_n
\) $ satisfies $ S_1^m (g_n) = \min \( \frac m n, 1 \) $, $ S_m^\infty
(g_n) = \max \( \frac{n-m+1}{n}, 0 \) $. According to the definitions of
Sect.~1, $ (f_n)_{n=1}^\infty $ is stable, but $ (g_n)_{n=1}^\infty $
is sensitive. We see that the definitions are not
tree-invariant. A straightforward way to tree-invariance is used in
the following definition of ``tree stability'' and ``tree
sensitivity''. From now on, stability and sensitivity of Sect.~1 will
be called ``cube stability'' and ``cube sensitivity''.

\begin{Def}
(a) A sequence $ (f_n)_{n=1}^\infty $ of functions $ f_n : \{-1,+1\}^n
\to \C $ is \emph{tree stable,} if there exists a sequence of tree
automorphisms $ A_n : \{-1,+1\}^n \to \{-1,+1\}^n $ such that the
sequence $ \( f_n \circ A_n \)_{n=1}^\infty $ is cube stable.

(b) The sequence $ (f_n)_{n=1}^\infty $ is \emph{tree sensitive,} if $
\( f_n \circ A_n \)_{n=1}^\infty $ is cube sensitive for every
sequence $ (A_n) $ of tree automorphisms.
\end{Def}

The definition can be formulated in terms of $ f_n ( A_n (\tau') ) $
and $ f_n ( A_n (\tau'') ) $ where $ (\tau',\tau'') $ is a $ \rho
$-correlated pair of random points of the cube $ \{-1,+1\}^n
$. Equivalently, we may consider $ f_n (\tau') $ and $ f_n(\tau'') $
where $ \tau', \tau'' $ are such that for some $ A_n $, $ ( A_n \tau',
A_n \tau'' ) $ is a $ \rho $-correlated pair. That is,
\begin{align}
& \E \( \tau'_m \big| \tau'_1, \tau''_1, \dots, \tau'_{m-1},
\tau''_{m-1} \) = \E \( \tau''_m \big| \tau'_1, \tau''_1,
\dots, \tau'_{m-1}, \tau''_{m-1} \) = 0 \, , \label{immers} \\
& \E \( \tau'_m \tau''_m \big| \tau'_1, \tau''_1, \dots, \tau'_{m-1},
  \tau''_{m-1} \)
 = a (\tau'_1,\dots,\tau'_{m-1}) a (\tau''_1,\dots,\tau''_{m-1}) \rho
\, , \label{factor}
\end{align}
where $ a : \cup_{m=1}^n \{-1,+1\}^{m-1} \to \{-1,+1\} $.
On the other hand, consider an arbitrary $ \{-1,+1\}^n \times
\{-1,+1\}^n $-valued random variable $ (\tau', \tau'') $ satisfying
\eqref{immers} (which implies that each one of $ \tau', \tau'' $ is
uniform on $ \{-1,+1\}^n $), but maybe not \eqref{factor}, and define
\begin{equation}
\rhomax (\tau', \tau'') = \max_{m=1,\dots,n} \max \big| \E \( \tau'_m
\tau''_m \big| \tau'_1, \tau''_1, \dots, \tau'_{m-1}, \tau''_{m-1} \)
\big| \, ,
\end{equation}
where the internal maximum is taken over all possible values of $ (
\tau'_1, \tau''_1, \dots, \linebreak[0] \tau'_{m-1}, \tau''_{m-1} ) $.
The joint distribution of $ \tau' $ and $ \tau'' $ is a probability
measure $ \mu $ on $ \{-1,+1\}^n \times \{-1,+1\}^n $, and we denote $
\rhomax ( \tau', \tau'' ) $ by $ \rhomax (\mu) $. Given $ f,g :
\{-1,+1\}^n \to \C $, we denote $ \E \overline{f(\tau')} g(\tau'') $
by $ \langle f | \mu | g \rangle $.

\begin{Def}
A sequence $ \( f_n \)_{n=1}^\infty $ of functions $ f_n : \{-1,+1\}^n
\to \C $, satisfying $ 0 < \liminf_{n\to\infty} \| f_n \| \le
\limsup_{n\to\infty} \| f_n \| < \infty $, is \emph{cosy,} if for any
$ \eps > 0 $ there is a sequence $ (\mu_n)_{n=1}^\infty $, $ \mu_n $
being a probability measure on $ \{-1,+1\}^n \times \{-1,+1\}^n $,
such that $ \limsup_{n\to\infty} \rhomax (\mu_n) < 1 $ and $
\limsup_{n\to\infty} \( \| f_n \|^2 - \langle f_n | \mu_n | f_n
\rangle \) < \eps $.
\end{Def}

\begin{lemma}
Every tree stable sequence is cosy.
\end{lemma}

\textsc{Proof.}
Let $ (f_n) $ be tree stable. Take tree automorphisms $ A_n $ such
that $ (f_n \circ A_n) $ is cube stable. We have $ \E \(
\overline{ f_n(A_n(\tau'))} f_n(A_n(\tau'')) \) \xrightarrow[\rho\to1]{}
\E | f_n (\tau) |^2 $ uniformly in $ n $. Here $ \tau', \tau'' $ are $ \rho
$-correlated. The joint distribution $ \mu_n (\rho) $ of $ A_n (\tau')
$ and $ A_n (\tau'') $ satisfies $ \rhomax (\mu_n(\rho)) \le \rho $
due to \eqref{factor}. Also, $ \langle f_n | \mu_n | f_n \rangle
\xrightarrow[\rho\to1]{} \| f_n \|^2 $ uniformly in $ n $, which means
that $ \sup_n \( \| f_n \|^2 - \langle f_n | \mu_n | f_n \rangle \)
\to 0 $ for $ \rho \to 1 $. \qed

\sk

Is there a cosy but not tree stable sequence? We do not know. The
conditional correlation given by \eqref{factor} is not only $ \pm \rho
$, it is also factorizable (a function of $ \tau' $ times the
same function of $ \tau'' $), which seems to be much stronger than
just $ \rhomax (\mu) \le \rho $.

\section{Hypercontractivity}

Let $ (\tau', \tau'') $ be a $ \rho $-correlated pair of random points
of the cube $ \{-1,+1\}^n $. Then for every $ f,g : \{-1,+1\}^n \to \R
$
\begin{equation}\label{hyper}
\big| \E f(\tau') g(\tau'') \big|^{1+\rho} \le \( \E
|f(\tau')|^{1+\rho} \) \( \E |g(\tau'')|^{1+\rho} \) \, ,
\end{equation}
which is a discrete version of the celebrated hypercontractivity
theorem pioneered by Nelson (see \cite[Sect.~3]{Nel}).
For a proof, see \cite{Ba}; there, following Gross \cite{Gr}, the
inequality is proved for $ n = 1 $ (just two points, $ \{-1,+1\} $)
\cite[Prop.~1.5]{Ba}, which is enough due to tensorization \cite[Lemma
1.3]{Ba}. (See also \cite[Lemma 2.4]{BKS}.) The case of $ f,g $ taking
on two values $ 0 $ and $ 1 $ only is especially important:
\[
\Prob^{1+\rho} \( \tau' \in S' \,\&\, \tau'' \in S'' \) \le \Prob ( \tau' \in
S' ) \Prob ( \tau'' \in S'' ) = \frac{|S'|}{2^n} \cdot \frac{|S''|}{2^n}
\]
for any $ S', S'' \subset \{-1,+1\}^n $. Note that $ \rho = 0 $ means
independence,\footnote{%
Equality results from the inequality applied to complementary sets.}
while $ \rho = 1 $ is trivial: $ \Prob^2 (\ldots) \le \( \min (
\Prob(S'), \Prob(S'') ) \)^2 \le \Prob (S') \Prob (S'') $.

For a probability measure $ \mu $ on $ \{-1,+1\}^n \times \{-1,+1\}^n
$ we denote by $ \langle g | \mu | f \rangle $ the value $ \E \(
f(\tau') g(\tau'') \) $, where $ (\tau', \tau'') \sim \mu $. The
hypercontractivity \eqref{hyper} may be written as $ \big| \langle g |
\mu | f \rangle \big| \le \| f \|_{1+\rho} \| g \|_{1+\rho} $, where $
\mu = \mu(\rho) $ is the distribution of a $ \rho $-correlated
pair. The class of $ \mu $ that satisfy the inequality (for all $ f,g $) is
invariant under transformations of the form $ A \times B $, where $
A,B : \{-1,+1\}^n \to \{-1,+1\}^n $ are arbitrary invertible maps 
(since such maps preserve $ \| \cdot \|_{1+\rho} $). In
particular, all measures of the form (\ref{immers}--\ref{factor}) fit.

Can we generalize the statement for all $ \mu $ such that $ \rhomax
(\mu) \le \rho \, $? The approach of Gross, based on tensorization,
works on cubes (and other products), not trees. Fortunately, we have
another approach, found by Neveu \cite{Nev}, that works also on
trees.

\begin{lemma}\label{hc}
For every $ r \in [\frac12,1] $, $ x,y \in [0,1] $, and $ \rho \in [ - 
\frac{1-r}r, \frac{1-r}r ] $,
\[\begin{split}
& (1+\rho) (1-x)^r (1-y)^r + (1-\rho) (1-x)^r (1+y)^r + \\
&\qquad + (1-\rho) (1+x)^r (1-y)^r + (1+\rho) (1+x)^r (1+y)^r \le 4 \, .
\end{split}\]
\end{lemma}

\textsc{Proof.}
The left hand side is linear in $ \rho $ with the coefficient $ \(
(1+x)^r - (1-x)^r \) \( (1+y)^r - (1-y)^r \) \ge 0 $. Therefore, it
suffices to prove the inequality for $ \rho = \frac{1-r}r $, $ r \in
(\frac12, 1) $ (the cases $ r = \frac12 $ and $ r=1 $ follow by
continuity). Assume the contrary, then the continuous function $ f_r $
on $ [0,1] \times [0,1] $, defined by
\[\begin{split}
 f_r (x,y) = \frac1r (1-x)^r (1-y)^r + \frac{2r-1}r (1-x)^r (1+y)^r +\qquad&
  \\
{}+ \frac{2r-1}r (1+x)^r (1-y)^r + \frac1r (1+x)^r (1+y)^r &\, ,
\end{split}\]
has a global maximum $ f_r (x_0,y_0) > 4 $ for some $ r \in (\frac12,
1) $.  The
case $ x_0 = y_0 = 0 $ is excluded (since $ f_r (0,0) = 4 $). Also, $
x_0 \ne 1 $ (since $ \frac{\pd}{\pd x} \big|_{x=1-} f_r (x,y) = -\infty $)
and $ y_0 \ne 1 $. The new variables
\[
u = \frac{1+x}{1-x} \in [1,\infty) \, , \quad v = \frac{1+y}{1-y} \in
[1,\infty)
\]
will be useful. We have
\begin{equation}\label{3.2a}
\frac{1+x}{(1-x)^r (1-y)^r} \frac{\pd}{\pd x} f_r (x,y) = u^r v^r - u
- (2r-1) ( uv^r - u^r ) \, .
\end{equation}
For $ u = 1 $, $ v > 1$ the right hand side is $ 2 (1-r) (v^r-1) > 0 $;
therefore $ x_0 \ne 0 $ (since $(x_0,y_0)\ne (0,0)$), and similarly
$ y_0 \ne 0 $. So, $ (x_0,y_0) $ is an
interior point of $ [0,1] \times [0,1] $. The corresponding $ u_0, v_0
\in (1,\infty) $ satisfy $ u_0^r v_0^r - u_0 - (2r-1) (u_0 v_0^r -
u_0^r) = 0 $.  
By subtracting the same expression with $v_0$ switched with $u_0$, 
which also vanishes, we get
\[
v_0 - u_0 + (2r-1) ( u_0^r v_0 - u_0 v_0^r + u_0^r - v_0^r ) = 0 \, .
\]
Aiming to conclude that $ u_0 = v_0 $, consider the function $ u
\mapsto v_0 - u + (2r-1) ( u^r v_0 - u v_0^r + u^r - v_0^r ) $ on $
[1,\infty) $. It is concave,  and positive when $u=1$,
since
$
v_0 - 1 + (2r-1) (v_0 - 2v_0^r + 1 )
\geq v_0 - 1 + (2r-1) (v_0 - 2v_0 + 1 )
= (v_0-1)(2-2r)
$.
Therefore, the function cannot vanish more than once, and
$ u = v_0 $ is its unique root. So, $ u_0 = v_0 $.

It follows from \eqref{3.2a} that
\[
\frac{1+x}{(1-x)^{2r}} \cdot \frac12 \frac{\pd}{\pd x} f(x,x) = u^{2r}
- u - (2r-1) ( u^{r+1} - u^r ) \, ,
\]
therefore $ u_0 $ is a root of the equation $ u^{2r-1} - 1 - (2r-1) (
u^r - u^{r-1} ) = 0 $, different from the evident root $ u = 1
$. However, the function 
 $ u \mapsto u^{2r-1} - 1 - (2r-1) ( u^r - u^{r-1} ) $
is strictly monotone, since
\[
\frac1{2r-1} \frac{\pd}{\pd u} (\ldots) = u^{2r-2} - r u^{r-1} + (r-1)
u^{r-2} = u^{r-2} ( u^r - ru + r - 1 ) < 0
\]
due to the inequality $ u^r \le 1 + r (u-1) $ 
(which follows from concavity of $u^r$).
The contradiction completes the proof.\qed

\sk

\begin{theorem}\label{th3.3}
Let $ \rho \in [0,1] $, and $ \mu $ be a probability measure on $
\{-1,+1\}^n \times \{-1,+1\}^n $ such that\footnote{%
It is assumed that $ \mu $ satisfies \eqref{immers}; $ \rhomax $ was
defined only for such measures.}
$ \rhomax (\mu) \le \rho $. Then for every $ f,g : \{-1,+1\}^n \to \C
$
\[
\big| \langle g | \mu | f \rangle \big| \le \| f \|_{1+\rho} \| g
\|_{1+\rho} \, .
\]
\end{theorem}

\sk

\textsc{Proof.} Consider random points $ \tau', \tau'' $ of $
\{-1,+1\}^n $ such that $ (\tau',\tau'') \sim \mu $. We have two
(correlated) random processes $ \tau'_1, \dots, \tau'_n $ and $
\tau''_1, \dots, \tau''_n $. Consider the random variables
\[
M'_n = | f (\tau'_1,\dots,\tau'_n) |^{1/r} \, , \quad
M''_n = | g (\tau''_1,\dots,\tau''_n) |^{1/r} \, ,
\]
and the corresponding martingales
\[\begin{split}
M'_m &= \E \( M'_n \big| \tau'_1, \tau''_1, \dots, \tau'_m, \tau''_m
\) = \E \( M'_n \big| \tau'_1, \dots, \tau'_m \) \, , \\
M''_m &= \E \( M''_n \big| \tau'_1, \tau''_1, \dots, \tau'_m, \tau''_m
\) = \E \( M''_n \big| \tau''_1, \dots, \tau''_m \)
\end{split}\]
for $ m = 0,1,\dots,n $; the equalities for conditional expectations
follow from \eqref{immers}. For any $ m = 1,\dots,n $ and any values
of $ \tau'_1, \tau''_1, \dots, \tau'_{m-1}, \tau''_{m-1} $ consider
the conditional distribution of the pair $ (M'_m, M''_m) $. It is
concentrated at four points that can be written as\footnote{%
Of course, $ x $ and $ y $ depend on $ \tau'_1, \tau''_1, \dots,
\tau'_{m-1}, \tau''_{m-1} $.}
$ \( (1\pm x) M'_{m-1}, (1\pm y) M''_{m-1} \) $. The first ``$ \pm $''
depends only on $ \tau'_m $, the second on $ \tau''_m $ (given the
past); each of them is ``$-$'' or ``$+$'' equiprobably. They have some
correlation coefficient lying between $ (-\rho) $ and $ \rho $. Lemma
\ref{hc} gives
\[
4 \E \bigg( \bigg( \frac{M'_m}{M'_{m-1}} \frac{M''_m}{M''_{m-1}}
\bigg)^r \bigg| \dots \bigg) \le 4 \, ,
\]
where $ r = \frac1{1+\rho} $. Thus, $ \E \( (M'_m M''_m)^r \big| \dots
\) \le (M'_{m-1} M''_{m-1})^r $, which means that the process $ (M'_m
M''_m)^r $ is a supermartingale. Therefore, $ \E (M'_n M''_n)^r \le
(M'_0 M''_0)^r $, that is, 
\[\begin{split}
\E | f (\tau'_1,\dots,\tau'_n) g (\tau''_1,\dots,\tau''_n) |
&
\le \( \E | f (\tau'_1,\dots,\tau'_n)
|^{1/r} \)^r \cdot \( \E | g (\tau''_1,\dots,\tau''_n) |^{1/r} \)^r 
\\ & =
\| f \|_{1+\rho} \| g \|_{1+\rho} \, .
\end{split}\]
\qed

\section{The main result}

Return to the spider walk (Fig.~1b). It may be treated as a
complex-valued martingale $ Z $ (Fig.~2a), starting at the origin.
Take each step to have length $1$. The
set $ \Om_n^\spi $ of all $ n $-step trajectories of $ Z $ can be 
identified with the set
of leaves of a binary tree. The endpoint $ Z_n = Z_n (\om) $ of a
trajectory $ \om \in \Om_n^\spi $ is a complex-valued function on $
\Om_n^\spi $. Taking into account that $ \E | Z_n |^2 = n $, we ask
about tree stability of the sequence $ \( Z_n / \sqrt n
\)_{n=1}^\infty $.

\begin{figure}[hb]
\begin{center}
\includegraphics{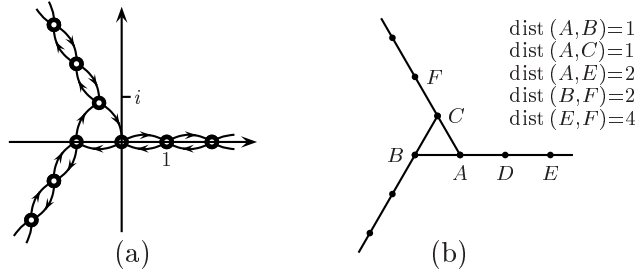}
\caption[]{(a) the spider walk as a complex-valued martingale;
(b) combinatorial distance.}
\end{center}
\end{figure}

\vspace{-0.5cm}

\begin{theorem}
The sequence $ \( Z_n / \sqrt n \)_{n=1}^\infty $ is non-cosy.
\end{theorem}

By Lemma 2.6 it follows that the sequence $ \( Z_n / \sqrt n
\)_{n=1}^\infty $ is not tree stable.
Recently, M.~Emery and J.~Warren found that some tree sensitive
sequences result naturally from their constructions.

In contrast to the spider walk, the simple walk (Fig.~1a) produces a
sequence $ \( (\tau_1+\dots+\tau_n)/\sqrt n \)_{n=1}^\infty $ that
evidently is cube stable, therefore tree stable, therefore cosy.

\begin{lemma}\label{lm}
(a) $ \limsup_{n\to\infty} \sqrt n \Prob ( Z_n = 0 ) < \infty $.

(b) $ \liminf_{n\to\infty} \( n^{-1/2} \sum_{k=1}^n \Prob(Z_k=0) \) >
0 $.
\end{lemma}

The proof is left to the reader.
Both (a) and (b) hold for each node of our graph, not just $ 0 $.
In fact, the limit exists, $
\lim_{n\to\infty} \( n^{1/2} \Prob(Z_n=0) \) = \frac12
\lim_{n\to\infty} \( n^{-1/2} \sum_{k=1}^n \Prob(Z_k=0) \) \in
(0,\infty) $, but we do not need it.

\medskip

\textsc{Proof of the theorem.}
Let $ \mu_n $ be a probability measure on $ \Om_n^\spi \times
\Om_n^\spi $ such that\footnote{%
It is assumed that $ \mu $ satisfies \eqref{immers}; $ \rhomax $ was
defined only for such measures.}
$ \rhomax (\mu) \le \rho $, $ \rho \in (0,1) $; we'll estimate $
\langle Z_n | \mu_n | Z_n \rangle $ from above in terms of $ \rho
$. We have two (correlated) copies $ \( Z'_k \)_{k=1}^n $, $ \( Z''_k
\)_{k=1}^n $ of the martingale $ \( Z_k \)_{k=1}^n $. Consider the
combinatorial distance (see Fig.~2b)
\[
D_k = \dist ( Z'_k, Z''_k ) \, .
\]
Conditionally, given the past $ ( Z'_1, Z''_1, \dots, Z'_{m-1},
Z''_{m-1} ) $, we have two \emph{equiprobable} values for $ Z'_m $,
and two \emph{equiprobable} values for $ Z''_m $; the two binary
choices are correlated, their correlation lying in $ [-\rho, \rho]
$. The four possible values for $ (Z'_m,Z''_m) $ lead usually to three
possible values $ D_{m-1} - 2 $, $ D_{m-1} $, $ D_{m-1} + 2 $ for $
D_m $, see Fig.~3a; their probabilities depend on the correlation, but the
(conditional) expectation of $ D_m $ is equal to $ D_{m-1} $
irrespective of the correlation. Sometimes, however, a different situation
appears, see Fig.~3b; here the conditional expectation of $ D_m $ is
equal to $ D_{m-1} + 1/2 $ rather than $ D_{m-1} $. That happens when $
Z''_{m-1} $ is situated at the beginning of a ray (any one of our
three rays) and $ Z'_{m-1} $ is on the same ray, outside the central
triangle $ \De $ ($ ABC $ on Fig.~2b). In that case\footnote{%
There is a symmetric case ($ Z'_{m-1} $ at the beginning\dots), but
we do not use it.}
we set $ L_{m-1} = 1 $, otherwise $ L_{m-1} = 0 $.
We do not care about the case when $ Z'_{m-1}, Z''_{m-1} $ are both on
$ \De $; this case may be neglected due to hypercontractivity, as
we'll see soon.  Also, the situation where $ Z'_{m-1}= Z''_{m-1} $ 
may occur, and then $\E\(D_m\big|D_{m-1}\)\ge D_{m-1}$.

\begin{figure}[ht]
\begin{center}
\includegraphics{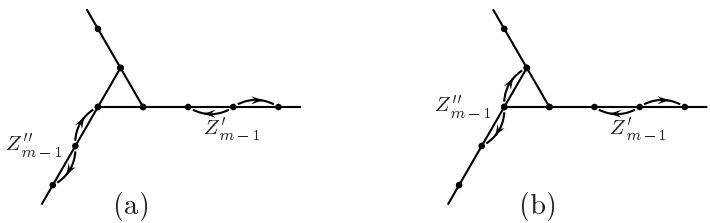}
\caption[]{(a)~the usual case, $ L=0 $: in the mean, $ D $ remains the same;
(b)~the case of $ L=1 $: in the mean, $ D $ increase by $ 1/2 $. More
cases exist, but $ D $ never decreases in the mean.}
\end{center}
\end{figure}

Theorem \ref{th3.3}, applied to appropriate indicators, gives $
\Prob^{1+\rho} \( Z'_k \in \De \,\&\, Z''_k \in \De \) \le \Prob \(
Z'_k \in \De \) \cdot \Prob \( Z''_k \in \De \) $, that is,
\[
\Prob \( Z'_k \in \De \,\&\, Z''_k \in \De \) \le \( \Prob ( Z_k \in
\De ) \)^{\frac2{1+\rho}}
\]
for all $ k = 0,\dots,n $. Combining it with Lemma \ref{lm} (a) we get
\begin{equation}\label{eqA}
\sum_{k=0}^n \Prob \( Z'_k \in \De \,\&\, Z''_k \in \De ) \le \eps_n
(\rho) \cdot \sqrt n
\end{equation}
for some $ \eps_n (\rho) $ such that $ \eps_n (\rho)
\xrightarrow[n\to\infty]{} 0 $ for every $ \rho \in (0,1) $, and $
\eps_n (\rho) $ does not depend on $ \mu $ as long as $ \rhomax (\mu)
\le \rho $.

Now we are in position to show that
\begin{equation}\label{eqB}
\sum_{k=0}^n \Prob \( L_k = 1 \) \ge c_0 \sqrt n
\end{equation}
for $ n \ge n_0 (\rho) $; here $ n_0 (\rho) $ and $ c_0 > 0 $ do not
depend on $ \mu $.  First, Lemma \ref{lm} (b) shows that $ \Prob (
Z''_k = 0 ) $ is large enough. Second, \eqref{eqA} shows that $ \Prob
( Z''_k = 0 \,\&\, Z'_k \notin \De ) $ is still large enough. The same
holds for $ \Prob ( Z''_k = 0 \,\&\, Z'_k \notin \De_{+2} ) $, where $
\De_{+2} $ is the (combinatorial) $ 2 $-neighborhood of $ \De $. Last,
given that $ Z''_k = 0 $ and $ Z'_k \notin \De_{+2} $, we have a
not-so-small (in fact, $ \ge 1/4 $) conditional probability that $ L_k
+ L_{k+1} + L_{k+2} > 0 $.  This proves \eqref{eqB}.

The process $ \( D_m - \frac12 \sum_{k=0}^{m-1} L_k \)_{m=0}^n $ is a
submartingale (that is, increases in the mean). Therefore, using
\eqref{eqB},
\[
\E D_n \ge \frac12 \sum_{k=0}^{n-1} \Prob ( L_k = 1 ) \ge \frac12 c_0 \sqrt n
\]
for $ n \ge n_0 (\rho) $. Note that $ D_n = \dist ( Z'_n, Z''_n ) \le
C_1 | Z'_n - Z''_n | $ for some absolute constant $ C_1 $. We have
\[
\( \E | Z'_n - Z''_n |^2 \)^{1/2} \ge \E | Z'_n - Z''_n | \ge C_1^{-1}
\E D_n \ge \frac12 C_1^{-1} c_0 \sqrt n
\]
and
\[
\| Z_n \|^2 - \langle Z_n | \mu_n | Z_n \rangle = \frac12 \E | Z'_n -
Z''_n |^2 \ge \frac14 C_1^{-2} c_0^2 n
\]
for $ n \ge n_0 (\rho) $; so,
\[
\limsup_{n\to\infty} \bigg( \bigg\| \frac{Z_n}{\sqrt n} \bigg\|^2 -
\bigg\langle \frac{Z_n}{\sqrt n} \bigg| \mu_n \bigg| \frac{Z_n}{\sqrt
n} \bigg\rangle \bigg) \ge \frac{c_0^2}{4 C_1^2}
\]
irrespective of $ \rho $, which means non-cosiness.\qed

\section{Connections to continuous models}

Theorem 4.1 (non-cosiness) is a discrete counterpart of
\cite[Th.~4.13]{Triple}. A continuous complex-valued martingale $ Z(t)
$ considered there, so-called Walsh's Brownian motion, is the limit of
our $ \( Z_{nt}/\sqrt n\) $ when $ n\to\infty $. The constants $ c_0 $ and
$ C_1 $ used in the proof of Theorem 4.1 can be improved (in fact, made
optimal) by using explicit calculations for Walsh's Brownian motion.
Cosiness for the simple walk is a discrete counterpart of \cite[Lemma
2.5]{Triple}.

Theorem 3.3 (hypercontractivity on trees) is a discrete counterpart
of \cite[Lemma 6.5]{Triple}. However, our use of hypercontractivity
when proving non-cosiness follows \cite[pp.~278--280]{BEKSY}. It is
possible to estimate $ \Prob ( Z'_k \in \De \,\&\, Z''_k \in \De ) $
without hypercontractivity, following \cite{EY} or
\cite[Sect.~4]{Triple}.

Cosiness, defined in Def.~2.5, is a discrete counterpart of the notion
of cosiness introduced in \cite[Def.~2.4]{Triple}. Different variants of 
cosiness (called I-cosiness and D-cosiness) are investigated by \'Emery,
Schachermayer, and Beghdadi-Sakrani, see \cite{E} and references
therein. See also Warren \cite{War}.

Noise stability and noise sensitivity, introduced in \cite{BKS}, have
their continuous counterparts, see \cite[Sect.~2]{Tsi}. Stability
corresponds to white noises, sensitivity to black noises
\cite{TV,Tsir}. Intermediate cases (neither stable nor sensitive, see
\cite[end of Sect.~1.4]{BKS}) correspond to noises that are neither
white nor black \cite{Warr}.

\bigskip
\filbreak
\begingroup
{
\small
\begin{sc}
\parindent=0pt\baselineskip=12pt
\def\emailwww#1#2{\par\qquad {\tt #1}\par\qquad {\tt #2}\medskip}

The Weizmann Institute of Science,
Rehovot 76100, Israel
\emailwww{schramm@wisdom.weizmann.ac.il}
{http://www.wisdom.weizmann.ac.il/$\sim$schramm/}

School of Mathematics, Tel Aviv Univ., Tel Aviv
69978, Israel
\emailwww{tsirel@math.tau.ac.il}
{http://math.tau.ac.il/$\sim$tsirel/}
\end{sc}
}
\filbreak

\endgroup


\begin{thebibliography}{99}

\bibitem{Ba} D. Bakry, ``L'hypercontractivit\'e et son utilisation en th\'eorie des semigroups'',
{\em Lect. Notes Math} (Lectures on probability theory),
Springer, Berlin, 1581 (1994), 1--114.

\bibitem{BEKSY} M.T. Barlow, M. \'Emery, F.B. Knight, S. Song, M. Yor,
``Autour d'un th\'eor\`eme de Tsirelson sur des filtrations browniennes
et non browniennes'',
{\em Lect. Notes Math} (S\'eminaire de Probabilit\'es XXXII),
Springer, Berlin, 1686 (1998), 264--305.

\bibitem{BKS} I. Benjamini, G. Kalai, O. Schramm, ``Noise
sensitivity of Boolean functions and applications to percolation'',
math.PR/9811157.

\bibitem{E} M.~\'Emery, ``Remarks on an example studied by
A.~Vershik and M.~Smorodinsky'', Manuscript, 1998.

\bibitem{EY} M.~\'Emery, M.~Yor,
``Sur un th\'eor\`eme de Tsirelson relatif \`a des
 mouvements browniens corr\'el\'es et \`a la nullit\'e
 de certains temps locaux'',
{\em Lect.\ Notes Math.} (S\'eminaire de Probabilit\'es XXXII),
Springer, Berlin, 1686 (1998), 306--312.

\bibitem{Gr} L. Gross, ``Logarithmic Sobolev inequalities'',
Amer. J. Math. 97 (1976), 1061--1083.

\bibitem{Nel} E. Nelson, ``The free Markoff field'',
J. Funct. Anal. 12 (1973), 211--227.

\bibitem{Nev} J. Neveu, ``Sur l'esp\'erance conditionelle par rapport
\`a un mouvement brownien''.
Ann. Inst. H. Poincar\'e 12 (1976), 105--109.

\bibitem{Triple} B. Tsirelson, ``Triple points: from non-Brownian
filtrations to harmonic measures,'' Geom. Funct. Anal. (GAFA) 7
(1997), 1096--1142.

\bibitem{Tsi} B.~Tsirelson,
``Unitary Brownian motions are linearizable'',
MSRI Preprint No.\ 1998-027, math.PR/9806112.

\bibitem{Tsir} B.~Tsirelson,
``Brownian coalescence as a black noise'',
manuscript in preparation.

\bibitem{TV} B.S.~Tsirelson, A.M.~Vershik,
``Examples of nonlinear continuous tensor products of
measure spaces and non-Fock factorizations'',
Reviews in Mathematical Physics 10:1 (1998), 81--145.

\bibitem{War} J.~Warren,
``On the joining of sticky Brownian motion'',
{\em Lect.\ Notes Math.} (S\'eminaire de Probabilit\'es XXXIII),
Springer, Berlin (to appear).

\bibitem{Warr} J.~Warren,
``The noise made by a Poisson snake'',
Manuscript, Univ.\ de Pierre et Marie Curie, Paris, Nov.\ 1998.

\end{thebibliography}
\end{document}